\documentclass[12pt,a4paper]{amsart}
\usepackage{amsmath,amssymb,enumerate,amsthm}
\usepackage{graphicx}
\usepackage{fancyhdr}

\theoremstyle{plain}

\input{tcilatex}

\begin{document}

\begin{center}
\bigskip {\Huge Classification trees in a box extent lattice\bigskip }

Laura Veres

Institute of Mathematics, University of Miskolc

3515 Miskolc-Egyetemv\'{a}ros, Hungary

matvlaura@gmail.com
\end{center}

\noindent \textit{\noindent Abstract. }

In this paper we show that during an elementary extension of a context each
of the classification trees of the newly created box extent lattice can be
obtained by the modification of the classification trees of the box extent
lattice of the original, smaller context. We construct also an algorithm
which, starting from a classification tree of the box extent lattice of the
smaller context $(H,M,I\cap H\times M),$ gives a classification tree of the
extended context $(G,M,I)$ which contains the new elements inserted. The
effectiveness of the method is that it ensures that there is enough to know
the original context, the classification tree of the box extent lattice and
its box extents, we do not need a new box extension of the extended context
mesh elements (except for one, which is the new element box extension).



\bigskip

\begin{center}
\bigskip {\Large 1.1 Preliminaries: box lattice, extent lattice\medskip }
\end{center}

A context (see [2]) is a triple $(G,M,I)$ where $G$ and $M$ are sets and $%
I\subseteq $ $G\times M$ is a binary relation. The elements of $G$ and $M$
are called\emph{\ objects} and respectively \emph{attributes} of the
context. The relation $gIm$ means that the object $g$ has the attribute $m$.
A small context can be easily represented by a cross table, i.e., by a
rectangular table, the rows of which are headed by the object names and the
columns are headed by the attribute names. A cross in the intersection of
the row $g$ and the column $m$, means that the object $g$ has the attribute $%
m$.

For all sets $A\subseteq G$ and $B\subseteq M$ we define

$\qquad A^{\prime }=\{m\in M\mid g$ $I$ $m$ for all $g\in A\}$,

$\qquad B^{\prime }=\{g\in G\mid g$ $I$ $m$ for all $m\in B\}$.

A\emph{\ concept} of the context\textbf{\ } $(G,M,I)$ is a pair $(A,B)$, in
which $A^{\prime }=B$ and $B^{\prime }=A$, and $A\subseteq G,$ $B\subseteq M$%
. $\tciLaplace (G,M,I)$ denotes the set of all concepts of the context $%
(G,M,I)$.

$\tciLaplace (G,M,I)$ can be endowed with the structure of a complete
lattice defining the join and meet of concepts as follows:

$\qquad \qquad \underset{i\in I}{\wedge }(A_{i},B_{i})=\left( \underset{i\in
I}{\cup }A_{i},\left( \underset{i\in I}{\cap }B_{i}\right) ^{\prime \prime
}\right) $

$\qquad \qquad \underset{i\in I}{\vee }(A_{i},B_{i})=\left( \left( \underset{%
i\in I}{\cap }A_{i}\right) ^{\prime \prime },\underset{i\in I}{\cup }%
B_{i}\right) $

The lattice $\left( \tciLaplace (G,M,I),\wedge ,\vee \right) $ will be
called the \emph{concept lattice} of the context $(G,M,I)$.

An \emph{extent partition} of a formal context $(G,M,I)$ is a partition of $%
G $, all classes of which are concept extents. Clearly, the trivial
partition $\{G\}$ is an extent partition. Note that, since the intersection
of extents always yields an extent, the common refinements of extent
partitions are still extent partitions. Therefore, the extent partitions of $%
(G,M,I)$ form a complete $\wedge $-subsemilattice of the partition lattice
of $G$, and thus a complete lattice which will be denoted with Ext$(G,M,I)$.
In particular, there is always a finest extent partition of the context
denoted with $\pi _{\square }.$

The zero element of a complete lattice $L$ and all elements that are
contained in some classification system of $L$ are called the \emph{box
elements} of $L$. The set of all box elements of $L$ is denoted by $Box(L)$
and $(Box(L),\leq )$ is a poset obtained by restricting the partial order of 
$L$ to the box elements. If every nonzero element of $L$ is a join of atoms
of $L$, then $L$ is called an \emph{atomistic lattice}.

In [10] the following was shown: If $L$ is a complete lattice in which every
element is a join of some completely join-irreducible elements, then $%
(Box(L),\leq )$ is a complete atomistic lattice.

In [3] the \emph{box extents} of a context $(G,M,I)$ where characterized if $%
E$ belongs to some extent partition of $(G,M,I)$ or $E=\varnothing ^{\prime
\prime }$. The set of all box extents of $L$ is denoted with $\mathcal{B}%
(G,M,I).$ The box extents are ordered by inclusion and in [3] was proved
that $\mathcal{B}(G,M,I)$ is a complete atomistic lattice which is a $\wedge 
$-subsemilattice of the concept lattice. (Note that if $\varnothing ^{\prime
\prime }\neq \varnothing $, then $\left\{ G\right\} $ is the only extent
partition of the context $(G,M,I)$). \noindent Observe, that each object $%
g\in G$\ is contained in a smallest box extent denoted with $g^{\square
\square }$, which is a class containing $g$\ of the finest extent partition $%
\pi _{\square }$\ of the context $(G,M,I)$.\ In other words\ $E$ is a box
extent iff $g\in E\Longrightarrow g^{\square \square }\in E$ or $g\notin
E\Longleftrightarrow g^{\square \square }\cap E=\varnothing .$\noindent
\noindent

\noindent\ \ \ \ \ \ In [4] was introduced the notion of CD-independent sets
in an arbitrary poset, we will define it in a lattice as follows:

\ Let $L$ be a bounded lattice. A set $X$ $\subseteq $\textit{\ }$\mathit{L}$
is called \emph{CD-independent}, if for any $x,y\in X$ either $x\leq y$ or $%
y\leq x$ or $x\wedge y=0$ (any elements of $X$ are comparable or disjoint).
Maximal CD-independent sets (with respect to $\subset $) are called \emph{%
CD-bases}.

Let $L$ be a lattice with the smallest element 0. A set $O=\{a_{i}|\,i\in
I\} $, $I\neq \varnothing $ of nonzero elements of $L$ is called a \emph{%
disjoint set} or \emph{orthogonal system}, if $a_{i}\wedge a_{j}=0$, $i\neq
j $. $O$ is a \emph{maximal orthogonal system}, if there is no other
orthogonal system $O^{\prime }$ of $L$ containing $O$ as a proper subset.
Notice that if $S_{1}=\{a_{i}|i\in I\}$, $I\neq \varnothing $ and $%
S_{2}=\{b_{j}|j\in J\}$, $J\neq \varnothing $ are two orthogonal systems,
then $S_{1}\leq S_{2}$, if for each $i\in I$ there exists $j(i)\in J$ such
that $a_{i}\leq b_{j(i)}$. We denoted with Ort$(L)$ the set of all
orthogonal systems, $(Ort(L),\leq )$ is evidently a poset, moreover a
lattice.

\noindent \textbf{Theorem 1.1:} ([4],[12])\noindent \emph{\ Let }$L$ \emph{%
be a lattice, and} $T$\emph{\ a CD-base in }$L$\emph{. }

\noindent (i) \emph{Then there exists a chain} $\mathcal{C}=\left\{
S_{\lambda }\mid \lambda \in \Lambda \right\} $\emph{\ in }Ort$(L)$\emph{,
such that }$T=\underset{\lambda \in \Lambda }{\cup }S_{\lambda }$\emph{. \ }

\noindent (ii) \emph{Any CD-base T\ is the union of disjoint sets
(orthogonal systems)} \emph{belonging to a maximal chain.}

\noindent

\noindent \bigskip

\begin{center}
{\LARGE 2. Classification trees and complete classification trees}
\end{center}

\bigskip

\noindent \qquad The classification trees are used in many fields: lattice
theory, data mining, group technology problems. Classification trees are
used for clustering the objects by their attributes, and they appear in some
clustering problems originated in Group Technology. Several independent
definitions are for classification trees in lattice theory, we will define
it as a special chain and also we will show some construction theorems for
them.

First we will relive the notion of classification trees and the relation
between classification trees, disjoint (orthogonal) systems and
CD-independent sets [11]. \medskip

Let $L$ be a lattice with the greatest element $1\in L$ and $T\subseteq
L\smallsetminus \{0\}$, $T\neq \emptyset $ a subset of it. $T$ \emph{is a
directed tree}, if $\left[ t\right) \cap T$ is a chain for each $t\in T$ and 
$1\in T$ .

\noindent Now strengthening the above condition we obtain the notion of a
classification tree:

\noindent \textbf{Definition 2.1.} $T$ is a \emph{classification tree}, if $%
\left[ x\right) \cap T$ is a chain (nonempty) for each $x\in L\setminus
\{0\} $. $H$ is a \emph{maximal classification tree}, if there is no other
classification tree which contains it as a proper subset.

Observe that in the definition of a classification tree the condition $\left[
t\right) \cap T$ is a chain is replaced with a stronger one: $\left[
x\right) \cap T$ is a chain (nonempty). As a consequence, we get that any
classification tree is also a directed tree but in addition for any $x,y\in
X $ either $x\wedge y=0$ or $x,y$ are comparable. Hence any classification
tree is CD-independent, and any maximal classification tree is a CD-base.
Moreover in [4,10] are proved the following assertions, which are also true
for classification trees and CD-independent sets:

\noindent \textbf{Proposition 2.2.} Let $L$ be a bounded lattice and $%
T\subseteq L$ a nonempty subset of it. Then the following are equal:

(i) $T$ is a maximal classification tree in $L$;

(ii) $T\cup \{0\}$ is a CD-base of $L$.

\noindent \textit{Remark 2.3: (i)} If $L$ is an atomistic lattice, $A(L)$
the set of its all atoms then any CD-base as well as any maximal
classification trees of $L$ contains $A(L)$.

\noindent (ii) If $T$ is a classification tree and we add atoms to it, it
remains a classification tree. (If $S\subseteq A(P)$, then $T\cup S$ is a
classification tree.) [11]

The notion of context and classification trees appears in Group Technology
problems. This engineering discipline exploits similarities between
technological objects and divides them into relatively homogenous groups,
extent partitions (classes) in order to optimize manufacturing processes.
The concept of classification tree appear also in the Group Technology. In
Group Technology by the term "classification tree" we refer to a tree made
of extents belonging to the context. This trees have the property that any
maximal antichain selected from the tree is a cover of the set $G$ with
extents.

Obvious that basis, the following definition results:

\noindent \textbf{Definition 2.4. }A classification tree\textbf{\ }$\mathcal{%
T}$ $\subseteq $ Ext$(G,M,I)$ \ is called a \emph{complete classification
tree} if for any maximal antichains $\left\{ E_{i}|\text{ }i\in I\right\}
\subseteq \mathcal{T}$ we have $\underset{i\in I}{\cup }E_{i}=G$. $\mathcal{T%
}$ is a \emph{maximal complete classification tree}, if there is no other
complete classification tree which contains it as a proper subset.

It's obvious that all elements of such a complete classification tree are
box extents, ie $\mathcal{T}$ is a classification tree in the lattice of box
extents.

\noindent \textbf{Proposition 2.5.}\emph{\ Let }$(G,M,I)$ \emph{be a finite
context,} $Ext(G,M,I)$ \emph{the \emph{extent }lattice of the context and} $%
\mathcal{B}(G,M,I)$ \emph{the lattice of the box extents of the context.
Then the followings are equal:}

\noindent (i) $\mathcal{T}\subseteq \;$Ext$(G,M,I)$ \emph{is a complete
classification tree; }

\noindent (ii) $\mathcal{T}$ \emph{is a classification tree} $\mathcal{B}%
(G,M,I)$\emph{\ and each maximal antichain of it }$\left\{ E_{i}|\text{ }%
i\in I\right\} \subseteq \mathcal{T}$ \ \emph{is a complete orthogonal
system in }$\mathcal{B}(G,M,I)$\emph{. }

\noindent \textit{Proof. }(i)$\Rightarrow $(ii) First we will prove that $%
\mathcal{T}$ is a classification tree in the box extent lattice $\mathcal{B}%
(G,M,I)$. Since $\mathcal{B}(G,M,I)$ is a subsemilattice of Ext$(G,M,I)$ and
it contains the smallest element ($\emptyset $) of Ext$(G,M,I)$, it is
enough to prove that each element $A\in \mathcal{T}$ \ is a box extent. As $%
\mathcal{T}$ is finite, then $A$ must be an element of a maximal antichain $%
\left\{ E_{i}|\text{ }i\in I\right\} \subseteq \mathcal{T}$, ie. $A=E_{k}$
for some $k\in I$. Since the elements of an antichain are incomparable, for
any two elements $E_{i},E_{j}\in \mathcal{T}$, $i,j\in I$ we have $E_{i}\cap
E_{j}=\emptyset $. As $\mathcal{T}$ is a complete classification tree in Ext$%
(G,M,I)$, we get $\underset{i\in I}{\cup }E_{i}=G$. Then $\left\{ E_{i}|%
\text{ }i\in I\right\} $ is an extent partition of $G$. Hence $A=E_{k}$ is a
box extent. From the above follows that any maximal antichain $E_{i}$, $i\in
I$ is a part of an orthogonal system, so it is an orthogonal system itself.
Obviously, this orthogonal system can not be extended with a new element $%
E_{0}\in \mathcal{B}(G,M,I)$, $E_{0}\neq \emptyset $. If we extended it,
then $E_{i}\cap E_{0}=\emptyset $ for each $i\in I$, from which follows that 
$E_{0}=G\cap E_{0}=$ $\left( \underset{i\in I}{\cup }E_{i}\right) \cap
E_{0}=\emptyset ,$ contradiction.

\noindent (ii)$\Rightarrow $(i) If $\mathcal{T}$ is a classification tree in
the lattice $\mathcal{B}(G,M,I)$, then it is also a classification tree in
Ext$(G,M,I)$. Assume that (ii) is satisfied, then it is enough to prove that
for every maximal antichain $\left\{ E_{i}|\text{ }i\in I\right\} \subseteq 
\mathcal{T}$ we have $\underset{i\in I}{\cup }E_{i}=G$. By our hypothesis $%
\left\{ E_{i}|\text{ }i\in I\right\} $ is a complete orthogonal system. Let $%
g\in G$ and $g^{\square \square }$ the box extents of it. As $\left\{ E_{i}|%
\text{ }i\in I\right\} \cup g^{\square \square }$ is not an orthogonal
system, then $g^{\square \square }\cap E_{k}\neq \emptyset $ for some $k\in
I $. Since $g^{\square \square }\,$\ is an atom in the lattice of box
extents, implies that $g^{\square \square }\subseteq E_{k}$. Hence $%
G\subseteq $ $\underset{i\in I}{\cup }E_{i}$. As $\underset{i\in I}{\cup }%
E_{i}\subseteq G$, $\underset{i\in I}{\cup }E_{i}=G.\hfill \square $

\noindent\textbf{\ Corollary 2.6. }$\mathcal{T}\subseteq \mathcal{B}(G,M,I)$ 
\emph{is a maximal classification tree in }$\mathcal{B}(G,M,I)$ \emph{if and
only if} $\mathcal{T}$\emph{\ is a maximal complete classification tree is
the extent lattice} Ext$(G,M,I)$\emph{.}\medskip

\noindent \textit{Proof. }Let $\mathcal{T}$ be a maximal classification tree
in $\mathcal{B}(G,M,I)$. Then $\mathcal{T}$ contains all atoms of the
lattice $\mathcal{B}(G,M,I).$ Let $S$ be a maximal antichain in $\mathcal{T}$%
.\ As $\mathcal{T}$ \ is a CD-independent set, then $S$ is an orthogonal
system. We have to prove that $S$ is a complete orthogonal system.

We denote with $A(\mathcal{B}(G,M,I))$ the set of all atoms of the lattice $%
\mathcal{B}(G,M,I).$ It is easy to see that $A(\mathcal{B}(G,M,I))$ is a
complete orthogonal system. In [4] was proved that an orthogonal system $O$
is complete if $A(L)\leqslant O$. (where $A(L)$ is the set of all atoms of
the lattice $L$ and $O$ is an orthogonal system). Thus it's enough to prove
that $A(\mathcal{B}(G,M,I)\leqslant S$. Indeed, if this inequation is not
satisfied then there exists an atom $a$ \ in $\mathcal{B}(G,M,I)$ which is
not smaller then any element of $S$. If $a\in \mathcal{T}$, then $S\cup
\{a\} $ is an antichain in $\mathcal{T}$, which is a contradiction. Thus we
have $A(\mathcal{B}(G,M,I))\leqslant S$ and $S$ is a maximal orthogonal
system. By the Proposition 2.5 $\mathcal{T}$ is a complete classification
tree in the lattice Ext$(G,M,I)$. Now, let $\mathcal{F}\subseteq \;$Ext$%
(G,M,I)$ be a complete classification tree which contains $\mathcal{T}$.
Hence, using Proposition 2.5, $\mathcal{F}$ is also a classification tree in
the lattice $\mathcal{B}(G,M,I)$ and $\mathcal{T}$ \ by definition is a
maximal classification tree in $\mathcal{B}(G,M,I)$, so we obtain $\mathcal{F%
}=\mathcal{T}$. Thus $\mathcal{T}$ is a maximal complete classification tree
in Ext$(G,M,I)$.

Conversely, assume that $\mathcal{T}\subseteq $Ext$(G,M,I)$ is a maximal
complete classification tree. Hence using Proposition 2.5 $\mathcal{T}$ is
also a classification tree in the box extent lattice $\mathcal{B}(G,M,I)$
Thus it is a subset of a maximal classification tree $\mathcal{M}$ of$%
\mathcal{\ B}(G,M,I).$ Therefore, $\mathcal{M}$ is a complete classification
tree in the lattice Ext$(G,M,I)$ and $\mathcal{T}\subseteq \mathcal{M}$ is
maximal, so we obtain $\mathcal{T}=\mathcal{M}$.

Finally we have that $\mathcal{T}$ is a maximal classification tree in the
lattice $\mathcal{B}(G,M,I)$.$\hfill \square $

\bigskip

\begin{center}
\bigskip \newpage

{\LARGE 3. Classification trees in the lattice of box extents}
\end{center}

In what follows we will consider a formal context $K=(G,M,I)$ with a
concrete and fixed technical meaning, where $G$ and $M$ are finite and
nonempty sets, $G$ denotes a fixed set of technical objects, $M$ denotes a
fixed set of some possible, technically relevant properties, and for any $%
g\in G$ and any $m\in M,$ $gIm$ means that the object (part) $g$ has the
property $m$. Additionally we suppose that the context does not contain rows
or columns filled with only zeros. A full zero column means that the
corresponding property not held by any of the parts, so it is irrelevant; a
full zero row corresponds to a part g possessing none of the properties from
our list.

\noindent\ \ \ \ \ \ \ Let$\ K=(G,M,I)$ be a context and\emph{\ }$%
K_{H}=(H,M,I\cap H\times M)$ a subcontext of it, where $H\subseteq G$ and
there exists a $z\in G$ such that $H=$\emph{\ }$G\smallsetminus \left\{
z\right\} .$ In this section we use the results from the article [3], in
which is the authors studied what how will change the box extents after a
one-object extension of the context. In [3] was shown that the intersection
of box extents is also a box extent. In [3] was also proved that any extent
partition of the subcontext is also an extent partition of $K$ and the box
extents of the subcontext are also box extents of $K$ and the following
propositions:

\bigskip

\noindent \textbf{Proposition} \textbf{3.1. } \emph{If }$\pi =\left\{ A_{k}%
\text{ }|\text{ }k\in K\right\} $\emph{\ is an extent partition of}$\ K$%
\emph{\ then }%
\begin{equation*}
\pi _{H}=\left\{ A_{k}\cap H\text{ }|\text{ }k\in K\right\} \setminus
\left\{ \emptyset \right\}
\end{equation*}%
\emph{\noindent is an extent partition of the subcontext }$(H,M,I\cap
H\times M).$\emph{\ (If }$\pi $\emph{\ is the finest extent partition of } $%
K $ \emph{then }$\pi _{H}$\ \emph{is called the restriction of the} \emph{%
extent partition }$\pi $\emph{\ and conversely.} $\pi _{H}$\emph{\ is not
necessarily the finest extent partition of }$H$\emph{). }

\noindent \textbf{Corollary 3.2.} [3] \emph{If }$E$\emph{\ is a box extent
of }$(G,M,I)$\emph{\ then }$E\cap H$\emph{\ is a box extent of }$(H,M,I\cap
H\times M).\medskip $

\noindent \textbf{Corollary 3.3.} [10] $\mathcal{B}(G,M,I)$\emph{\ is a
complete atomistic lattice. The atomic box extents are the classes of the
finest extent partition }$\pi _{\square }$\emph{.\medskip }

\noindent \smallskip \textbf{Proposition} \textbf{3.4. }[3] \emph{If }$E$%
\emph{\ is a box extent of }$(G,M,I)$\emph{\ and }$H=G\smallsetminus \left\{
z\right\} $\emph{\ for some }$\emph{z}\in G$ \emph{then}

\noindent (1)\emph{\ }$E$\emph{\ is a box extent of }$(H,M,I\cap H\times M)$%
\emph{\ with }$E\cap \ z^{\square \square }=\emptyset $\emph{\ or}

\noindent (2)\emph{\ }$E\smallsetminus \left\{ z\right\} $\emph{\ is a box
extent of }$(H,M,I\cap H\times M)$\emph{.}\medskip

\noindent \textbf{Proposition 3.5. }[3]\textbf{\ }\emph{If }$H=$\emph{\ }$%
G\smallsetminus \left\{ z\right\} $ \emph{then} $A$ \emph{is a class of
finest extent partition of }$\pi _{\square }$ \emph{of} $(G,M,I)$ \emph{if
and only if }

\noindent (1)\emph{\ either} $A=z^{\square \square }$\emph{, or}

\noindent (2)\emph{\ } $A$ \emph{is a class, disjoint from }$z^{\square
\square }$\emph{,} \emph{of the finest extent partition of the subcontext }$%
(H,M,I\cap H\times M)$\emph{.}\bigskip

\noindent \textbf{Theorem 3.6. }[3] \emph{Let }$(G,M,I)$ \emph{be a context, 
}$E$ \emph{a box extent of the subcontext }$(H,M,I\cap H\times M)$ \emph{%
with }$H=$\emph{\ }$G\smallsetminus \left\{ z\right\} .$ \emph{Then}

\noindent (i)\emph{\ }$E$ \emph{is a box extent of} $(G,M,I)$ \emph{if and
only if } $z^{\square \square }\cap E^{\prime \prime }=\emptyset ;$

\noindent (ii)\emph{\ }$E^{\ast }=E\cup \left\{ z\right\} $ \emph{is a box
extent of }$(G,M,I)$ \emph{if and only if } $z^{\square \square }\setminus
\left\{ z\right\} \subseteq E$ \emph{and} $\left( E\cup \left\{ z\right\}
\right) ^{\prime \prime }=E\cup \left\{ z\right\} $.\bigskip

\noindent \textit{Remark:} In other words, in view of the theorem below, we
have two possibilities:

1) E is also box extent in the new context iff $z^{\square \square }\cap
E^{\prime \prime }=\emptyset $;

2) or $E\cup \left\{ z\right\} $ is a box extent in the new context iff $%
z^{\square \square }\setminus \left\{ z\right\} \subseteq E$ \emph{and} $%
\left( E\cup \left\{ z\right\} \right) ^{\prime \prime }=E\cup \left\{
z\right\} .\medskip $

\noindent $E^{\ast }=$\noindent $\QATOPD\{ . {E,\text{ if }z^{\square
\square }\cap E^{\prime \prime }=\emptyset \text{ }}{E\cup \left\{ z\right\} 
\text{, else}}\medskip $

In the following we show that during an elementary extension of the context
each of the classification trees of the newly created box extent lattice can
be obtained by the modification of the classification trees of the box
extent lattice of the original, smaller context.

\noindent \textbf{Proposition 3.7. }\emph{Let} $K_{H}=(H,M,I\cap H\times M)$ 
\emph{be} \emph{a subcontext of} $K=(G,M,I)$ \emph{and} $\mathcal{T}$ \emph{%
a classification tree in the box extent lattice }$\mathcal{B}(K).$\emph{\
Then the set}

\begin{center}
$\mathcal{T}_{H}=\{E\cap H\mid E\in\mathcal{T\}}$
\end{center}

\noindent \emph{is a classification tree in the box lattice} $\mathcal{B}%
(K_{H})$\emph{.}\medskip

\noindent \textit{Proof.} In [3] was proved that if $E$\ is a box extent of
the context $(G,M,I)$\ then $E\cap H$\ is the box extent of the subcontext $%
(H,M,I\cap H\times M).$\emph{\ \ }Thus for any $E\in \mathcal{T}$, $E\cap H$
is a box extent of the context $K_{H}$, so $\mathcal{T}_{H}\subseteq 
\mathcal{B}(K_{H})$. Since $G\in \mathcal{T}$, then $H=G\cap H\in \mathcal{T}%
_{H}$ and also $H$ is the largest element of the box extent lattice $%
\mathcal{B}(K_{H}).$ We have to prove that $\mathcal{T}_{H}$ is a
classification tree. It is enough to prove that $\mathcal{T}_{H}$ is a
CD-independent set of the lattice $\mathcal{B}(K_{H}).$ Let $E_{1}\cap H$
and $E_{2}\cap H$ be two incomparable elements of $\mathcal{T}_{H}$. Then $%
E_{1}\cap H\neq \emptyset $ and $E_{2}\cap H\neq \emptyset $ and $%
E_{1},E_{2}\in \mathcal{T}$ are also incomparable. Since $\mathcal{T}$ is
classification tree and also a CD-independent set, we have $E_{1}\cap
E_{2}=\emptyset .$ Then $(E_{1}\cap H)\cap (E_{2}\cap H)=E_{1}\cap E_{2}\cap
H=\emptyset $. Thus we proved that $\mathcal{T}_{H}$ is also a
CD-independent set in the lattice $\mathcal{B}(K_{H})$. Therefore $\mathcal{T%
}_{H}$ is a classification tree in $\mathcal{B}(K_{H}).\hfill \square $%
\bigskip

Further we show how to construct a classification tree after a one-object
extension of a context:

\noindent \textbf{Theorem 3.8.}\emph{\ Let }$K_{H}=(H,M,I\cap H\times M)$ 
\emph{be subcontext of a finite context }$(G,$ $M,I)$\emph{\ such that }$%
H=G\setminus \left\{ z\right\} ,$ $z^{\square \square }\neq \{z\}.$ \emph{Let%
} $\mathcal{T}$ \emph{be a classification tree in the box extent lattice} $%
\mathcal{B}(K_{H}).$\emph{\ Then:}

\noindent (i) $\mathcal{T}^{(1)}=\left\{ E\in \mathcal{T}|\text{ }E\in 
\mathcal{B}(G,M,I)\right\} $\emph{\ is an order ideal} \emph{in} $\mathcal{T}
$\emph{\ and} \newline
$\mathcal{T}^{(2)}=\left\{ E\in \mathcal{T}|\text{ }E\cup \left\{ z\right\}
\in \mathcal{B}(G,M,I)\right\} $ \emph{a finite chain in} $\mathcal{T}$\emph{%
\ and }$\mathcal{T}^{(1)}\cap \mathcal{T}^{(2)}=\emptyset $\emph{;}

\noindent (ii) $\mathcal{T}^{\ast }=\mathcal{T}^{(1)}\cup \left\{ E\cup
\left\{ z\right\} |\text{ }E\in \mathcal{T}^{(2)}\right\} $\emph{\ is a
classification tree in the lattice }$\mathcal{B}(G,M,I)$\emph{.}

\noindent (iii) \emph{If } $\mathcal{T}$ \emph{contains all the atoms of the
box extent lattice} $\mathcal{B}(K_{H})$\emph{, then }$\mathcal{T}^{\ast
}\cup \{z^{\square \square }\}$ \emph{is a classification tree in }$\mathcal{%
B}(G,M,I)$\emph{, which contains all the atoms of it.}

\noindent \textit{Proof. }(i) Let $E\in \mathcal{T}^{(1)}$ and $F\subseteq E,%
\mathcal{\ }F\in \mathcal{T}$. Since $E$ is a box extent of $(G,M,I)$ in
view of Theorem 3.6 we get $z^{\square \square }\cap E^{\prime \prime
}=\emptyset $. Since $F^{\prime \prime }\subseteq E^{\prime \prime }$ we
have $z^{\square \square }\cap F^{\prime \prime }=\emptyset $, which in view
of Theorem 3.6 means that $F\in \mathcal{T}^{(1)}$. Thus $\mathcal{T}^{(1)}$
is an order ideal in $\mathcal{T}$. We have to prove now that $\mathcal{T}%
^{(2)\text{ }}$is a finite chain in $\mathcal{T}$. We take the set $\mathcal{%
C}=\left\{ E\in \mathcal{T}|\text{ }z^{\square \square }\setminus \left\{
z\right\} \subseteq E\right\} $. Obviously, $H\in \mathcal{C}$, then $%
\mathcal{C}\neq \emptyset .$ As $z^{\square \square }\setminus \left\{
z\right\} \neq \emptyset $ we get $E\cap z^{\square \square }\neq \emptyset
. $ Thus $\mathcal{T}^{(1)}$ and $\mathcal{C}$ have no common elements. Take 
$E\in \mathcal{C}$ and $F\in \mathcal{T}$. Since$\ E\subseteq F$, we have $%
z^{\square \square }\setminus \left\{ z\right\} \subseteq F$ , so $F\in 
\mathcal{C}$. Therefore $\mathcal{C}$ is an order filter of $\mathcal{H}$
and since $H$ is finite and $\mathcal{C\subseteq }$ $\mathcal{H}$ is lower
bounded by it's minimal elements. We show that $\mathcal{C}$ has only one
minimal element$\ E_{1}\in \mathcal{C}$. Assume that $E_{2}\in \mathcal{C}$
is a minimal element in $\mathcal{C}$ and $E_{1}\neq E_{2}$. Since $%
E_{1},E_{2}\in \mathcal{T}$ and $E_{1},E_{2}$ are incomparable, we have $%
E_{1}\cap E_{2}=\emptyset $, which is a contradiction because $z^{\square
\square }\setminus \left\{ z\right\} \subseteq E_{1}\cap E_{2}$ and $%
z^{\square \square }\setminus \left\{ z\right\} $ are nonempty by
hypothesis. Thus $E_{1}$ is the smallest element of $\mathcal{C}$ and $%
\mathcal{C}$ is equal to $\left[ E_{1}\right) \cap \mathcal{T}$. Since $%
\mathcal{T}$ is a classification tree $\mathcal{C}$ must be a chain.
Consider now the set $\mathcal{T}^{(2)}$. By Proposition 3.5 $\mathcal{T}%
^{(2)}$ $\subseteq $ $\mathcal{C}$, then $\mathcal{T}^{(2)}$ is a finite
chain and $\mathcal{T}^{(1)}\cap \mathcal{T}^{(2)}=\emptyset $.

\noindent (ii) As $\mathcal{T}^{(1)}\subseteq \mathcal{T}$ we have that $%
\mathcal{T}^{(1)}$ is a CD-independent set. Observe that $\left\{ E\cup
\left\{ z\right\} \mid E\in \mathcal{T}^{(2)}\right\} $ is a chain in $%
\mathcal{B}(G,M,I)$. Indeed, let $E_{1}\cup \left\{ z\right\} $ and $%
E_{2}\cup \left\{ z\right\} $ be two elements of this set. Since $\mathcal{T}%
^{(2)}$ is a chain, then we have $E_{1}\subseteq E_{2}$ or $E_{2}\subseteq
E_{1}$ and $E_{1}\cup \left\{ z\right\} \subseteq E_{2}\cup \left\{
z\right\} $, or conversely $E_{2}\cup \left\{ z\right\} \subseteq E_{1}\cup
\left\{ z\right\} $. We show that the set $\mathcal{T}^{(1)}\cup \left\{
E\cup \left\{ z\right\} \mid E\in \mathcal{T}^{(2)}\right\} $ is also
CD-independent.

\noindent Take $E_{1}\in \mathcal{T}^{(1)}$ and $E_{2}\in \mathcal{T}^{(2)}$%
. As $E_{1},E_{2}\in \mathcal{T}$ \ $\mathcal{T}$ is a classification tree,
we have the following three cases: $E_{1}\subseteq E_{2}$ or $E_{2}\subseteq
E_{1}$ or $E_{1}\cap E_{2}=\emptyset $. We have to show that also the sets $%
E_{1}$ and $E_{2}\cup \left\{ z\right\} $ are either comparable or disjoint.

Clearly, in the first case $E_{1}\subseteq E_{2}\cup \left\{ z\right\} $.

In the second case $E_{2}\subseteq E_{1}\ $implies $E_{2}\in \mathcal{T}%
^{(1)}$ \ and $\mathcal{T}^{(1)}$ is an order ideal in $\mathcal{T}$.
However, this is impossible, because $E_{2}\in \mathcal{T}^{(2)}$\ and $%
\mathcal{T}^{(1)}\cap \mathcal{T}^{(2)}=\emptyset $.

Let us consider now the case $E_{1}\cap E_{2}=\emptyset $. Since $E_{1}\in 
\mathcal{T}^{(1)}$ we have $E_{1}^{\prime \prime }\cap z^{\square \square
}=\emptyset \ $and this results also $E_{1}\cap \left\{ z\right\} =\emptyset
.$Therefore, we obtain $E_{1}\cap \left( E_{2}\cup \left\{ z\right\} \right)
=\emptyset $.

Thus we have proved that $\mathcal{T}^{\ast }=\mathcal{T}^{(1)}\cup \left\{
E\cup \left\{ z\right\} \mid E\in \mathcal{T}^{(2)}\right\} $ is
CD-independent.

Since $G=H\cup \left\{ z\right\} $ and $H\in \mathcal{T}$ ( and $G\in 
\mathcal{B}(G,M,I)$), we obtain $G\in \mathcal{T}^{\ast }$. As $G$ is the
greatest element of the lattice $\mathcal{B}(G,M,I)$, in view of Remark 2.3 $%
\mathcal{T}^{\ast }$ is a classification tree in $\mathcal{B}(G,M,I)$.

\noindent (iii) Observe that because $z^{\square \square }$ is an atom in
the lattice $\mathcal{B}(G,M,I)$, if we add it to the classification tree $%
\mathcal{T}^{\ast }$, then in view of Proposition 2.2 $\mathcal{T}^{\ast
}\cup \{z^{\square \square }\}$ remains a classification tree in $\mathcal{B}%
(G,M,I)$.

Finally assume that $\mathcal{T}$ contains all the atoms of the lattice $%
\mathcal{B}(K_{H}).$\emph{\ }We have to show that the classification tree $%
\mathcal{T}^{\ast }\cup \{z^{\square \square }\}$ contains all the atoms of
the lattice $\mathcal{B}(G,M,I)$. Evidently, it contains the atom $%
z^{\square \square }$ also. On the other hand the atoms of $\mathcal{B}%
(G,M,I)$ are blocks of the finest extent partition $\pi _{\square }$ of $%
(G,M,I).$ Then in view of \ Proposition 3.5 all the other atoms of $\mathcal{%
B}(G,M,I)$ which are different from $z^{\square \square }$ are also atoms of
the lattice $\mathcal{B}(K_{H})$ , and belongs to $\mathcal{T}$. Therefore,
this atoms belongs to $\mathcal{T}^{(1)}$ by the construction of $\mathcal{T}%
^{(1)}$. Since $\mathcal{T}^{(1)}\subseteq \mathcal{T}^{\ast }\cup
\{z^{\square \square }\}$, then $\mathcal{T}^{\ast }\cup \{z^{\square
\square }\}$ contains all the atoms of $\mathcal{B}(G,M,I)$.$\hfill \square $

\medskip \noindent \textbf{Theorem 3.9. }\emph{Let }$K_{H}=(H,M,I\cap
H\times M)$ \emph{be a subcontext of the finite context }$K=(G,$ $M,I)$\emph{%
\ such that }$H=G\setminus \left\{ z\right\} ,$ $z^{\square \square }\neq
\{z\}$. \emph{Then:}

\noindent (i) \emph{For each classification tree} $\mathcal{T}_{G}\subseteq 
\mathcal{B}(G,M,I)$ \emph{there exists a classification tree} $\mathcal{T}%
_{H}\subseteq \mathcal{B}(K_{H})$ \emph{such that the equality} $\mathcal{T}%
_{G}=\mathcal{T}_{H}^{\ast }$ \emph{is satisfied.}

\noindent (ii) \emph{If} $\mathcal{T}_{G}\subseteq \mathcal{B}(G,M,I)$\emph{%
\ is a maximal classification tree, then} \emph{in }$\mathcal{B}(K_{H})$ 
\emph{exists also a maximal classification tree }$\mathcal{M}$\emph{\ such
that} $\mathcal{T}_{G}=\mathcal{M}^{\ast }$ \emph{.}\medskip

\noindent \textit{Proof.} (i) In view of Proposition 3.7 $\mathcal{T}%
_{H}=\{E\cap H\mid E\in \mathcal{T}_{G}\mathcal{\}}$ is a classification
tree in $\mathcal{B}(K_{H})$. We show that the classification tree

$\mathcal{T}_{H}^{\ast}=\mathcal{T}_{H}^{(1)}\cup\left\{ E\cup\left\{
z\right\} |\text{ }E\in\mathcal{T}_{H}^{(2)}\right\} =\{E\in\mathcal{T}%
_{H}\mid E\in\mathcal{B}(G,M,I)\}\cup$

$\left\{ E\cup\left\{ z\right\} |\text{ }E\in\mathcal{T}_{H}\text{, }%
E\cup\left\{ z\right\} \in\mathcal{B}(G,M,I)\right\} \medskip$

\noindent assigned to $\mathcal{T}_{H}$ by Theorem 3.8 is equal to $\mathcal{%
T}_{G}$. ($\mathcal{T}_{H}^{(1)}$, $\mathcal{T}_{H}^{(2)}$ was defined in
Theorem 3.8)

\noindent Let $E\in \mathcal{T}_{G}$ be arbitrary. Then $F=E\cap H\in 
\mathcal{T}_{H}$ by the definition of $\mathcal{T}_{H}$. We have only two
cases:

(1) if $E\subseteq H$, then $E=F$;

(2) if $E\nsubseteq H$, then $E=F\cup \{z\}$.

\noindent In the first case $F\in \mathcal{T}_{H}^{(1)}$ so we have $E=F\in 
\mathcal{T}_{H}^{(1)}\subseteq \mathcal{T}_{H}^{\ast }$.

\noindent In the second case as $F\cup \{z\}\in \mathcal{B}(G,M,I)$ we have $%
E=F\cup \{z\}\in \mathcal{T}_{H}^{\ast }$.

\noindent In both cases $\mathcal{T}_{G}\subseteq \mathcal{T}_{H}^{\ast }$
because $E\in \mathcal{T}_{H}^{\ast }$ .

\noindent Conversely, let $E\in \mathcal{T}_{H}^{\ast }$ be arbitrary. By
the construction of $\mathcal{T}_{H}^{\ast }$ we have two cases or $E\in 
\mathcal{T}_{H}^{(1)}\subseteq \mathcal{T}_{H}$, or $E=F\cup \{z\}$, $F\in 
\mathcal{T}_{H}^{(2)}\subseteq \mathcal{T}_{H}$.

\noindent Then by the definition of $\mathcal{T}_{H}$ exists an $A\in 
\mathcal{T}_{G}$ such that in the first case we have $A\cap H=E$, and $A\cap
H=F$ in the second case.

\noindent First we show that $z\in A$ is not possible.

\noindent Indeed, in view of Theorem 3.6, as $E\in \mathcal{B}(K_{H})$, $%
E\cap z^{\square \square }=E^{\prime \prime }\cap z^{\square \square
}=\emptyset $ . On the other hand from $z\in A$ we would get $z^{\square
\square }\subseteq A$ so $z^{\square \square }\setminus \{z\}\subseteq A\cap
H=E$. Combining with the first result this would imply $z^{\square \square
}\setminus \{z\}=\emptyset $,ie $z^{\square \square }=\{z\}$, contradiction.

\noindent Thus in the case of $A\cap H=E$ we have $z\notin A$. Therefore, $%
A\subseteq H$ and $E=A\cap H=A\in \mathcal{T}_{G}$.

\noindent In the case two $A\cap H=F$. Observe that $A\subseteq H$ would
imply $F=A\in \mathcal{B}(G,M,I)$ which means that $F\in \mathcal{T}%
_{H}^{(1)}$. Howover this is not possible because by the definition $F\in 
\mathcal{T}_{H}^{(2)}$ and $\mathcal{T}_{H}^{(2)}\cap \mathcal{T}%
_{H}^{(1)}=\emptyset $. Therefore $A\nsubseteq H$, which means $z\in A$.
Then $A=F\cup \{z\}=E$. Thus we get $E=A\in \mathcal{T}_{G}$.

\noindent Since in both cases we proved $E\in \mathcal{T}_{G}$ , it follows $%
\mathcal{T}_{H}^{\ast }=\mathcal{T}_{G}$.

\noindent (ii) Assume that $\mathcal{T}_{G}$ is a maximal classification
tree in the lattice $\mathcal{B}(G,M,I).$Since $\mathcal{T}_{G}=\mathcal{T}%
_{H}^{\ast }$, if $\mathcal{T}_{H}$ is a maximal classification tree in $%
\mathcal{B}(K_{H})$, we done. If $\mathcal{T}_{H}$ is not maximal, then in $%
\mathcal{B}(K_{H})$ exists a maximal classification tree $\mathcal{M}$ in $%
\mathcal{B}(K_{H})$, such that $\mathcal{T}_{H}\subseteq \mathcal{M}$. Then
we have $\mathcal{T}_{G}=\mathcal{T}_{H}^{\ast }\subseteq \mathcal{M}^{\ast
} $. \ In view of Theorem 3.8 $\mathcal{M}^{\ast }$ is also a classification
tree in $\mathcal{B}(G,M,I)$. As $\mathcal{T}_{G}$ is a maximal
classification tree by hypothesis, we get $\mathcal{T}_{G}=\mathcal{M}^{\ast
}=\mathcal{T}_{H}^{\ast }$.\hfill $\square $

\bigskip

\bigskip

\begin{center}
{\LARGE 4. An algorithm for \ tree-construction in the box extent lattice}
\end{center}

\bigskip

In chapter 3 of the article, we examined that the box extents of a context
are not changed with a one-object extension or reduction. We add an object
to the context which has attributes from the existing $M$ attribute set of
the context. We showed that during a one-object extension of the context
each of the classification trees of the newly created box extent lattice can
be obtained by the modification of the classification trees of the box
extent lattice of the original, smaller context.

Based on Theorem 3.8 we construct an algorithm which, starting from a
classification tree of the box extent lattice of the smaller context $%
(H,M,I\cap H\times M),$ gives a classification tree of the extended context $%
(G,M,I)$ which contains the new elements inserted. The effectiveness of this
method is that it ensures that there is enough to know the original context,
the classification tree of the box extent lattice and its box extents, we do
not need a new box extent of the extended context mesh elements (except for
one, which is the box extent of the new element).

For the construction of the classification trees of the box extent lattice
we will use the recursive algorithm ORTOFA presented in [12]. For the
algorithm ORTOFA we need the original context $K$ and the box extents which
are contained in the matrix $DS$. The next step is to determine $z^{\square
\square }$. For finding the extent partition $E=z^{\square \square }$ ,
containing the new inserted element we will use Algorithm 2 presented by K%
\"{o}rei A. in [3].

\noindent The function Lista\_Fa first searches for the elements of the
order ideal part of our classification tree, ie. we check if $z^{\square
\square }$ is smaller then any box extent stored in the matrix $L$. If this
condition is satisfied we put the element in the matrix S1. After that using
the KOBJ(KTUL) we verify that the elements of S1 are in fact box extents,
and order them with the function. As a last step we find the chain part of
our classification tree.

\noindent In the algorithm we used two functions:

BENNE (A, B) - verify that the matrix A of the box extents contains or not
the elements of B;

BERAK (A, V, m) - a procedure that inserts in A the vector V as a new m +
1-th row.

In the function KOBJ(KTUL) the function KTUL finds the common properties of
the object set, and the function KOBJ finds the common objects of the
attribute set and used this two functions one after the other we got the set 
$A^{\prime \prime }$ of the object set $A$.

We use the following matrices:

DS(mxn) contains the box extents,

L(mxn) contains the copy of the box extents,

S1(mx(n+1)) contains the elements above to $z^{\square \square }$,

S(mxn) contains the chain above to $z^{\square \square }$ and F(mxn) stores
the elements below.

\bigskip

1. $\ LISTA\_FA(DS)\qquad \qquad $

2. $L\longleftarrow DS\qquad \qquad \qquad \qquad $

3. $\ S1\longleftarrow \varnothing \qquad \qquad \qquad \qquad $

4. $\ S\longleftarrow \varnothing \qquad \qquad \qquad \qquad $

5. $F\longleftarrow \varnothing $ \qquad \qquad \qquad \qquad

6. $\ k\longleftarrow 0$

7. $\ for$ $i\longleftarrow 1$ $to$ $m$ \qquad \qquad \qquad \qquad

8. $\ \qquad do$ \ \ $benne\longleftarrow true$

9. $\qquad \qquad for$ $j\longleftarrow 1$ $to$ $n$

10.\qquad \qquad \qquad $do$ $if$ $z^{\square \square }[j]>L[i][j]$

11.\qquad \qquad \qquad \qquad $then$ $benne\longleftarrow false$

12.\qquad \qquad $if$ $benne$

13.\qquad \qquad \qquad $then$ $k++$

14.$\qquad \qquad \qquad \qquad BERAK(S1,L[i],k)\qquad $

15. $for$ $i\longleftarrow 1$ $to$ $k$

16. \qquad $do$ $S1[i][n+1]\longleftarrow 1$

17. $A\longleftarrow S1^{T}$

18. $A\longleftarrow KOBJ(KTUL(A,K),K)\qquad \qquad $*/here K means the
extended context

19. $A\longleftarrow A^{T}$

20. $l\longleftarrow 0$

21. $for$ $i\longleftarrow 1$ $to$ $k$

22.$\qquad do$ \ \ $dext\longleftarrow true$

23.$\qquad \qquad for$ $j\longleftarrow 1$ $to$ $(n+1)$

24. \qquad \qquad \qquad $do$ $if$ $A[i][j]\neq S1[i][j]$

25. \qquad \qquad \qquad \qquad \qquad $then$ $dext\longleftarrow false$

26. \qquad \qquad $if$ $dext$

27. \qquad \qquad \qquad $then$ $l++$

28. \qquad \qquad \qquad \qquad $BERAK(S,S1[i],l)$

29. $h\longleftarrow 0$

30. $for$ $i\longleftarrow 1$ $to$ $m$

31. \qquad $do$ $if$ $BENNE(L[i],$ $z^{\square \square })$

32. \qquad \qquad \qquad $then$ $h++$

33. \qquad \qquad \qquad \qquad $BERAK(F,L[i],h)$

34. $return$ $S,F$

\bigskip

$BERAK(A,V,m)$

1. $m++$

2. $for$ $t\leftarrow 1$ $to$ $n$

3. \qquad $do$ $A[m][t]\longleftarrow V[t]$

4. $return$ $A,m$\bigskip

$BENNE(A,B)$

1.$bent\leftarrow true$

2.$for$ $i\leftarrow 1$ $to$ $n$

3.$\qquad do$ $if$ $A[i]>B[i]$

4.\qquad \qquad \qquad $then$ $bent\leftarrow false$

5.$return$ $bent$

The effectiveness of the method is that it ensures that there is enough to
know the original context, the classification tree of the box extent lattice
and its box extents, we do not need a new box extents of the extended
context mesh elements (except for one, which is the new element box extent).

The main parts of the process are the dual matrix operations operating
cycles. The run-time of this algorithm is $O(n^{2})$ polynomial time, this
means the worst run time for a sufficiently large $n=max\{n,m,k\}$ value.
However, the method contains many one-line instruction (decisions, variable
value increase with 1) if the value of n is not too big, then the running
time is linear, $O(n)$ because the number of steps, the executable
instructions are $n$. In summary, the worst running time is second-order
polynomial time.\bigskip

\textbf{References}

[1] F\"{o}ldes S., Radeleczki S.: On interval decomposition lattices.
Discussiones Mathematicae, General Algebra and Applications 24 (2004), p.
95-114.

[2] Ganter B., Wille R.: Formal Concept Analysis, Mathematical Foundations.
Springer Verlag, Berlin-Heidelberg-New York (1999).

[3] Ganter, B., K\"{o}rei, A. and Radeleczki, S.: Extent partitions and
context extensions, Math. Slovaca, 63 (4) (2013), 693--706.

[4] E. K. Horv\'{a}th and S. Radeleczki, Notes on CD-independent subsets,
Acta Sci. Math. (Szeged), 78 (2012), 3-24.

[5] K\"{o}rei A. and Radeleczki S.: Box elements in a concept lattice.
Contribution to ICFCA 2006, Editors Bernhard Ganter, Leonard Kwuida, Verlag
Allgemeine Wissenschaft, Drezda 2006.

[6] MITROFANOV, S. P.: Scientific Principles of Group Technology, National
Lending Library for Science and Technology, Boston, 1966.

[7] Radeleczki S.: Classification systems and their lattice. Discussiones
Mathematicae, General Algebra and Applications 22 (2002), p. 167-181.

[8] Radeleczki S.: Fogalomh\'{a}l\'{o}k \'{e}s alkalmaz\'{a}suk a
csoporttechnol\'{o}gi\'{a}ban. Proceedings of International Computer
Science, microCAD, Miskolc, 1999, 3-8.

[9] Radeleczki S.: Classification systems and the decomposition of a lattice
into direct products. Mathematical Notes, Vol.1. No.1, 2000, 145-156.

[10] Radeleczki S.: Classification systems and their lattice. Discussiones
Mathematicae, General Algebra and Applicationes 22, 2002, 167-181.

[11] Veres L.: Generalization of classification trees for a poset. Creative
Mathematics and Informatics, Vol.17, 2008, 72-77;

[12] Veres L.: Classification trees and orthogonal systems. microCAD 2007
International Scientific Conference, Vol. G, Miskolc, 2007, 63-68;

[13] Wille R.: Subdirect decomposition of concept lattices. Algebra
Universalis 17, 1983, 275-287.

\end{document}